# Mean, Variance and Asymptotic Property for General Hypergeometric Distribution

Xing-gang Mao[a*] and Xiao-yan Xue[b]

[a]Department of Neurosurgery, Xijing Hospital, the Fourth Military Medical University, Xi'an, No. 17 Changle West Road, Xi'an, Shaanxi Province, China;

[b]Department of Pharmacology, School of Pharmacy, the Fourth Military Medical University, No. 17 Changle West Road, Xi'an, Shaanxi Province, China

*Corresponding authors. E-mail addresses: xgmao@fmmu.edu.cn (X.G. Mao)

Short title: General hypergeometric distribution

**Abstract**

*General hypergeometric distribution* (GHGD) definition: from a finite space $N$ containing $n$ elements, randomly select totally $T$ subsets $M_i$ (each contains $m_i$ elements, $1 \leq i \leq T$), what is the probability that exactly $x$ elements are overlapped exactly $t$ times or at least $t$ times ($x_t$ or $x_{\geq t}$)? The GHGD described the distribution of random variables $x_t$ and $x_{\geq t}$. In our previous results, we obtained the formulas of mathematical expectation and variance for special situations ($T \leq 7$), and not provided proofs. Here, we completed the exact formulas of mean and variance for $x_t$ and $x_{\geq t}$ for any situation, and provided strict mathematical proofs. In addition, we give the asymptotic property of the variables. When the mean approaches to 0, the variance fast approaches to the value of mean, and actually, their difference is a higher order infinitesimal of mean. Therefore, when the mean is small enough (<1), it can be used as a fairly accurate approximation of variance.

**Keywords**: overlapped elements; statistical distribution; group allocation of particles; finite population; mean; variance; mathematical expectation

## 1 Introduction

We will use italic upper-case letters to represent a set, and the corresponding italic lower-case letter to represent the cardinality of the set, that is number of elements in the set. Specifically, we use $T$ to represent the number of subsets.

Definition of *general hypergeometric distribution* (**GHGD**) has been described in our previous report[1]: from a finite space $N$ containing $n$ elements, randomly select totally $T$ subsets $M_i$ (each contains $m_i$ elements, $1 \leq i \leq T$; here $m_i = s_i$), **what is the probability that exactly $x$ elements are overlapped exactly $t$ times or at least $t$ times ($x_t$ or $x_{\geq t}$)?** By mathematical form, what is the statistical distribution of the random variables $x_t$ and $x_{\geq t}$? The probabilities are denoted as:

$$p(x_t = k \mid n, m_1, m_2 ... m_T) \quad \text{or} \quad p(x_t = k \mid n, M) \quad [1]$$

and

$$p(x_{\geq t} = k \mid n, m_1, m_2 ... m_T) \quad \text{or} \quad p(x_{\geq t} = k \mid n, M) \quad [2]$$

Here $M = \{m_1, m_2, ... m_T\}$.

When $T = 2$, the GHGD reduced to classical hypergeometric distribution (HGD). For very long time, the GHGD is not fully studied. This issue can also be described as the problem of group allocation of particles, for which the number of empty cells (corresponding to the number of 0-overlapped elements, $x_0$) has been studied [2]. However, in general, although efforts have been made, no more progress was made except for our previous work[1], mainly because of its difficulty even as a pure mathematical problem, as partially reflected by the challenge to develop tools to plot the Venn diagram [3-6] even for very small amount of sets (such as *4 - 6* sets). Actually, statistically dealing with finite populations are commonly difficult tasks [7]. There are attempts to establish methods to perform statistical inference, but also limited to 2 sets comparison [8]. Due to a lack of GHGD theory, only a few studies tried to get the significance of overlapped genes, but restricted to 2 sets overlaps, or replaced by pairwise comparisons in multiple subsets [9-12]. Another study developed a procedure to calculate multi-set intersections by using the R software package, *SuperExactTest*, but mainly based on enumerating method [13], and the distribution theory is not established. The work also discussed the difficulty of this problem caused mainly by huge calculations, further highlighting the significance to establish a distribution theory to avoid unpractical calculations.

Considering the importance and widespread application prospect of GHGD, here we further studied and completed the basic statistical theory of GHGD. The main achievement of the present study include: 1. We discovered exact formulas for essential statistics of the GHGD, including the mathematical expect (or mean) and variance for any $T$ and $t$ values. 2. We demonstrated the asymptotic properties of GHGD when the mean is very small. In detail, when the mean approaches to 0, the variance fast approaches to the value of the mean, and actually, their difference is a higher order infinitesimal of the mean. Therefore, when the mathematical expectation is small enough (<1), it can be used as a fairly accurate approximation of variance. As a pure mathematic-based statistical theory, it can be used in a wide range of situations when it is needed to consider the overlapped elements in multiple subsets. Detailed mathematical proofs were provided separately in the appendix for readability.

## 2 Results

### 2.1 Definition and key statistics for GHGD

In classical HGD, the parameters for the distribution include the number of the space $n$, number of the two subsets $m_1$, $m_2$. Then the distribution is denoted as $p(x=k \mid n, m_1, m_2)$. Here $k$ takes all possible values of number of overlapped elements (NOE) of the two subsets. However, in GHGD, there are totally $T$ ($T \geq 2$) subsets that have $m_1$, $m_2$ … $m_T$ elements, respectively. We denote the number of elements in these subsets as an integer set $M$. The elements in these subsets can be overlapped for $0$ to $T$ times. Therefore, the parameters for the GHGD are: number of the space $n$, numbers of elements in the selected subsets $m_1$, $m_2$ … $m_T$, and the feature of the overlap, $LO=t$ or $LO \geq t$. Hence, the random variable is defined as $x_t$ and $x_{\geq t}$. The probability of the event "$x_t = k$" or "$x_{\geq t} = k$" defines the probability density distribution of $x_t$ and $x_{\geq t}$, which is denoted as equation [1] or [2]. Here, $k$ takes all possible values of NOE of the $T$ subsets, that is between 0 and $Min(\sum_{i=1}^{T} m_i, N)$.

To simplify, if not confused, the distributions are denoted as $p(x_t = k)$ and $p(x_{\geq t} = k)$. Similarly, the number of all possible subsets meeting the condition "$x_t = k$" or "$x_{\geq t} = k$" is simplified as $C(x_t = k)$ and $C(x_{\geq t} = k)$. The mathematical expectation and variance are labeled as $E(x_t)$ or $E(x_{\geq t})$, and $Var(x_t)$ or $Var(x_{\geq t})$.

The total number of all possible $T$ subsets containing $m_i$ elements drawn from the space $N$ is

labeled as $C_n^M$, which is calculated by:

$$C_n^M = \prod_{i=1}^{T} C_n^{m_i}$$

Therefore, the probability is calculated by:

$$p(x_t = k) = \frac{C(x_t = k)}{C_n^M} \text{ or } p(x_{\geq t} = k) = \frac{C(x_{\geq t} = k)}{C_n^M}.$$

### 2.1.1 Basic formulations

First, we give the basic formulations. When $t = T$, the random variable $x_T$, that is the number of fully overlapped elements, has a mathematical expectation as follows:

$$E(x_T) = \frac{\prod_{i=1}^{T} m_i}{n^{T-1}} \quad [6]$$

Its variance is:

$$var(x_T) = \frac{\prod_{i=1}^{T} m_i}{n^{T-1}} \times \left(1 + \frac{\prod_{i=1}^{T}(m_i - 1)}{(n-1)^{T-1}} - \frac{\prod_{i=1}^{T} m_i}{n^{T-1}}\right)$$

$$= E(x_T)(1 + E(x_T \mid n-1, M-1) - E(x_T)) \quad [7]$$

Here, $E(x_T, n\text{-}1, M\text{-}1)$ is the mathematical expression of distribution $p(x_T = k \mid n-1, m_1 - 1, m_2 - 1, \ldots m_T - 1)$. The variance is the 2rd central moments about the mean. Actually any central moments (*v*th) about the mean can be calculated (detailed mathematical proof of these formulas are provided in **Appendix 1**).

### 2.1.2 Mathematical expectation for $p(x_t = k)$ and $p(x_{\geq t} = k)$

Based on the basic formulas, the mathematical expectation of $p(x_t=k)$ for any $T$ and $t$ value, is calculated by the following formulas.

$$E(x_t) = \sum_{M_t \subset \{M_t\}} E(x_T \mid n, M(M_t)) = \sum_{M_t \subset \{M_t\}} E(x_T \mid n, M_t \cup \overline{M_t})$$

$$= \sum_{m_t \subset M_t} \frac{\prod_{i=1}^{t} m_{t_i} \prod_{i=1}^{d}(n - \overline{m_{t_i}})}{n^{T-1}} \quad [8]$$

Here $d = T\text{-}t$, $M_t$ is one of the *t*-subset (that is a subset containing *t* element) of $M$, and $\{M_t\}$ is the set of all possible *t*-subsets of $M$. As a critical definition, $\overline{M_t}$ is the **complementary set** of $M_t$ in

$M$, that is, $\overline{M_t} = \{n - \overline{m}_{t_1}, n - \overline{m}_{t_2}, \ldots n - \overline{m}_{t_d}\}$, here $\overline{m}_{t_i}$ ($1 \leq i \leq d$) is the elements that belong to the $M$ but not to the $M_t$, and $M(M_t)$ is the union set of Mt and $\overline{M_t}$. For example, take $M = \{m_1, m_2, m_3, m_4, m_5\}$, $t = 3$, let $M_t = \{m_1, m_2, m_4\}$ is one of the 3-subsets of $M$, then $\overline{M_t} = \{n-m_3, n-m_5\}$, and $M(M_t) = \{m_1, m_2, m_4, n-m_3, n-m_5\}$. Obviously, each $E(x_T|n, M(M_t))$ in the formula [8] can be calculated with the basic formula [6]. This is because the variable $x_t$ is the sum of variables $(x_t|n,M_t)$, that is the number of $t$-overlapped elements in all $t$-subsets in $M$:

$$x_t = \sum_{M_t \subset \{M_t\}} (x_t \mid n, M_t) \qquad [9]$$

Therefore $E(x_t)$ can be calculated by the sum of mathematical expectations of these variables. Detailed proofs and explanations are provided in **Appendix 2**.

When all $m_i$ take the same value $m$, then we have

$$E(x_t) = C_T^t \times \frac{(n-m)^{T-t} \times m^t}{n^{T-1}} \qquad [10]$$

Next, for the situation of $OL \geq t$ we have:

$$E(x_{\geq t}) = E(x_t) + E(x_{\geq (t+1)}) = \sum_{i=t}^{T} E(x_i) \qquad [11]$$

### 2.1.3 Variances for $p(x_t = k)$ and $p(x_{\geq t} = k)$

According to the equation [9], we further obtained the variances of the variables $x_t$ and $x_{\geq t}$. That is:

$$E(x_t^2) = E\left(\left(\sum_{M_t \subset \{M_t\}} (x_t \mid n, M_t)\right)^2\right) = \sum_{M_{t1}, M_{t2} \subset \{M_t\}} E((x_{t1} \mid n, M_{t1})(E(x_{t2} \mid n, M_{t2})) \qquad [12]$$

$$= \sum_{M_{t1}, M_{t2} \subset \{M_t\}} E((x_{T1} \mid n, M(M_{t1}))) (E(x_{T2} \mid n, M(M_{t2})))$$

Obviously, when $M_{t1}=M_{t2}$, the item in the summary is the second moments $E(x_T^2)$. The key problem is to calculate the item when $M_{t1} \neq M_{t2}$. Here we have the result:

$$E(y = x_{Mt} \times x_{Ms}) = E(y = (x_t \mid n, M_t) \times (x_s \mid n, M_s)) = \frac{\prod_{m_i \subset M_{noSame}} m_i(n - m_i) \prod_{m_j \subset M_{same}} m_j(m_j - 1)}{n^T (n-1)^{T-1}} \qquad [13]$$

Here $x_{Mt}$ and $x_{Ms}$ are the variables of number of $t$ or $s$-overlapped elements in the $t$ or $s$-subset $M_t$ and $M_s$ ($t$ and $s$ can be equal or not). Here is another critical definition: $M_{noSame}$ and $M_{same}$. $M_{same}$ is the intersection set of $M(M_t)$ and $M(M_s)$, and $M_{noSame}$ is the set of elements that are belong to the

$M(M_t)$ but not to the $M_{same}$. For example, take $M = \{m_1, m_2, m_3, m_4, m_5\}$, let $M_t = \{m_1, m_2, m_3\}$, $M_s = \{m_1, m_2, m_4\}$, we have $M_{same} = \{m_1, m_2, n-m_5\}$ and $M_{noSame} = \{m_3, n-m_4\}$.

When $t=s$, and $M_t$ takes all possible $t$-subsets of $M$, then the sum equation [12] can be calculated to get $E(x_t^2)$. In addition, when $t$ and $s$ takes all possible values from $t$ to $T$, then we obtained the formulas of $E(x_{\geq t}^2)$.

Based on this result, we obtained the formulas of $E(x_t^2)$ and $E(x_{\geq t}^2)$, and therefore we can calculate the variances by using $Var(x^2) = E(x^2) - E^2(x)$. Detailed proof of [13] was also provided in **Appendix 2**. Importantly, the correctness of the formulas was also validated by the programmable and exhaustive enumeration algorithms.

## 2.2 Asymptotic property of GHGD

We found that when $E(x)$ is small enough, the difference between $E(x)$ and $Var(x)$ is very small, especially when $E(x)<1$. Actually, the difference is a higher order infinitesimal of $E(x)$ when $E(x)$ approaches to 0 ($E(x)\to 0$) (**Appendix 3**). In addition, we have proved that, $E(x) > Var(x)$ (**Appendix 3**). Therefore, even for larger $E(x)$ values (such as <100), $E(x)$ is an upper estimation of $Var(x)$. Importantly, when $E(x)$ is small enough (<1), we can use $Var(x) \approx E(x)$.

Let $\lambda = 1 - \mu$, then the Chebyshev's inequality can be transformed to (considering $x > 1$): $P(x \geq 1) \leq \dfrac{\sigma^2}{(1-\mu)^2}$, which is transformed as (let $\sigma^2 = Var(x) \approx E(x)$):

$$P(x \geq 1) \leq \dfrac{E(x)}{(1-E(x))^2} \qquad [14]$$

More specifically, when $E(X) \ll 1$, it can be further simplified as $p(x \geq 1) \leq E(x)$, that is the $E(x)$ value can be used to approximately estimate the event $p(x \geq 1)$. Lastly, we can estimate the $E(x)$ value so that the $P(x \geq 1) \leq \dfrac{E(x)}{(1-E(x))^2}$ reaches to a significant level $\alpha$. Specifically, when $\alpha = 0.05$, the maximal value of $E(x)$ is about 0.04554 ($\approx 0.05$). A more accurate estimation of is based on the $3\sigma$ rule from the Vysochanskii-Petunin inequality for unimodal distribution[14, 15]: the probability is less than 5% that $x$ is more than three standard deviations from the mean of the population: $P(|x-\mu|>3\sigma) \leq 0.05$. Take $\sigma^2 = E(x)$, we have the maximal value of $E(x)$ is about 0.09167. Therefore, if we get a result of $E(x) \leq 0.09167$, we can tell that the probability $p(x \geq 1) < 0.05$.

## 3 Applications of GHGD

The GHGD can be used to infer a overlapped subsets with statistical significance, which has been described in our previous report[1]. Here, we will give a novel potential application of the GHGD, which is very important for modern research of high-order interactions among multiple variables beyond two variables. For two variables, their interactions can be calculated by correlation or mutual information based analysis. For multiple variables, each has a series of value, then a threshold can be set to get a subset of elements which has are larger or smaller than the threshold. Then GHGD can be used to examine whether the subsets of all variables derived from the threshold has overlapped elements with statistical significance. Very interestingly, this kind of statistical inference identified high-order interactions of several variables, which has no low-order interactions. For example, in multi-variable analysis, triplet variables can be identified which had significant elements overlapped in the 3 sets, but had no significant pair-wise interactions or associations. Detailed description of this application will be provided in our another report.

## 4 Appendix: Mathematical Proofs

### 4.1 Appendix 1: basic formulations: Distribution of fully overlapped elements $p(x_T = k)$

For the simplest situation, we first considered the distribution of $p(x_T=k)$, where $k$ represents the number of elements that are overlapped in all of the subsets. Because for $t = T$, the distribution of $x_{LO=T}$ and $x_{LO \geq T}$ is identical, it will be labeled as $x_T$. We discovered graceful formulas for essential statistics of $p(x_T=k)$, including mathematical expectation, variance, and high order moments. In this section, if not specified, the parameters for variable $x_t$ is $(n, m_1, m_2, \ldots m_T)$. In addition, $M$-$j$ is used to represent $m_1$-$j$, $m_2$-$j$, $\ldots m_T$-$j$. For example, $p(x_T=k) = p(x_T=k|n, m_1, m_2, \ldots m_T) = p(x_T=k|n, M$-$j)$.

First we give calculations for $C(x_T=k) = C(x_T=k|n, m_1, m_2, \ldots m_T) = C(x_T=k|n, M)$. $C(x_T=k)$ can be calculated with the following recursive formula (here $m_{min}$ is the minimum of all $m_i$ ($1 \leq i \leq T$); $C_m^n = \frac{m!}{n!(m-n)!}$ is the binomial coefficient. Note: $C_m^k = 0$ when $k > m$; $C_m^0 = 1$):

$$C(x_T = k) = C_n^k \times \prod_{i=1}^{T} C_{n-k}^{m_i - k} - \sum_{i=k+1}^{m_{min}} C_i^k \times C(x_T = i) \quad [A1]$$

When $T=2$, the distribution reduced to classical HGD, confirming the correctness of the formula.

**Lemma 1:** $C(x_T = i) = \frac{n}{i} C(x_T = i-1 | n-1, M-1)$.

Proof by mathematical induction:

First, it is true take $i=m_{min}$. this is because, according to equation [A1]:

$$C(x_T = m_{min} | n, M) = C_n^{m_{min}} \times \prod_{i=1}^{T} C_{n-m_{min}}^{m_i - m_{min}} = \frac{n}{m_{min}} C_{n-1}^{m_{min}-1} \times \prod_{i=1}^{T} C_{(n-1)-(m_{min}-1)}^{(m_i-1)-(m_{min}-1)}$$

$$= \frac{n}{m_{min}} C(x_T = m_{min} - 1 | n-1, M-1).$$

Next, if it is true when $i \geq k+1$, then for $i = k$ we have:

$$C(x_T = k | n, m_1, m_2, ... m_T) = C_n^k \times \prod_{i=1}^{T} C_{n-k}^{m_i-k} - \sum_{i=k+1}^{m_{min}} C_i^k \times C(x_T = i | n, m_1, m_2, ... m_T)$$

$$= \frac{n}{k} C_{n-1}^{k-1} \times \prod_{i=1}^{T} C_{n-1-(k-1)}^{m_i-1-(k-1)} - \sum_{i=k+1}^{m_{min}} \frac{i}{k} C_{i-1}^{k-1} \times \frac{n}{i} C(x_T = i-1 | n-1, m_1-1, m_2-1, ... m_T-1)$$

$$= \frac{n}{k} \left( C_{n-1}^{k-1} \times \prod_{i=1}^{T} C_{(n-1)-(k-1)}^{(m_i-1)-(k-1)} - \sum_{i=k+1}^{m_{min}} C_{i-1}^{k-1} C(x_T = i-1 | n-1, m_1-1, m_2-1, ... m_T-1) \right)$$

$$= \frac{n}{k} \left( C_{n-1}^{k-1} \times \prod_{i=1}^{T} C_{(n-1)-(k-1)}^{(m_i-1)-(k-1)} - \sum_{i=(k-1)+1}^{m_{min}-1} C_i^{k-1} C(x_T = i | n-1, m_1-1, m_2-1, ... m_T-1) \right)$$

$$= \frac{n}{k} C(x_T = k-1 | n-1, m_1-1, m_2-1, ... m_T-1) = \frac{n}{k} C(x_T = k-1 | n-1, M-1)$$

**Theorem 1**: the mathematical expectation of $p(x_T=k)$ is:

$$E(x_T) = \frac{\prod_{i=1}^{T} m_i}{n^{T-1}}$$

Proof:

$$E(x_T) = \sum_{k=0}^{m_{min}} k \times \frac{C_n^k \times \prod_{i=1}^{T} C_{n-k}^{m_i-k} - \sum_{i=k+1}^{m_{min}} C_i^k \times C(x_T = i | n, M)}{\prod_{i=1}^{T} C_n^{m_i}}$$

$$= \sum_{k=0}^{m_{min}} k \times \frac{\frac{n}{k} C_{n-1}^{k-1} \times \prod_{i=1}^{T} C_{n-k}^{m_i-k} - \sum_{i=k+1}^{m_{min}} \frac{i}{k} C_{i-1}^{k-1} \times \frac{n}{i} C(x_T = i-1 | n-1, M-1)}{\prod_{i=1}^{T} \frac{n}{m_i} \prod_{i=1}^{T} C_{n-1}^{m_i-1}}$$

$$= \frac{\prod_{i=1}^{T} m_i}{n^T} \sum_{k=0}^{m_{min}} \frac{k \times \frac{n}{k} C_{n-1}^{k-1} \times \prod_{i=1}^{T} C_{n-1-(k-1)}^{m_i-1-(k-1)} - \sum_{i=k+1}^{m_{min}} k \times \frac{n}{k} C_{i-1}^{k-1} \times C(x_T = i-1 | n-1, M-1)}{\prod_{i=0}^{T-1} C_{n-1}^{m_i-1}}$$

$$= \frac{\prod_{i=0}^{T-1} m_i}{n^{T-1}} \sum_{k=0}^{m_{min}} \frac{C_{n-1}^{k-1} \times \prod_{i=1}^{T} C_{n-1-(k-1)}^{m_i-1-(k-1)} - \sum_{i=k+1}^{m_{min}} C_{i-1}^{k-1} \times C(x_T = i-1 | n-1, M-1)}{\prod_{i=0}^{T-1} C_{N-1}^{M[i]-1}}$$

The sum formula is the sum of another distribution $p(x_T=k|n-1, M-1)$, therefore its value is 1, and the result is $\dfrac{\prod_{i=1}^{T} m_i}{n^{T-1}}$.

Notably, when all M[i] is the same value M, then $E_{m_i=m_0}(x_T) = \dfrac{m_0^T}{n^{T-1}}$.

**Theorem 2**: Any moments (vth) of $p(x_T)$ can be calculated by the recursion formula (considering that $E(x_T^0)=1$):

$$E(x_T^v) = E(x_T) \times \sum_{i=0}^{v-1} C_{v-1}^i E(x_T^i \mid n-1, M-1)$$

Proof:

First, when $v=0$, $E(x_T^0)=1$.

If it is true for $j$th ($j \leq v-1$) moments, then:

$$E(x_T^v) = \sum_{k=0}^{m_{\min}} k^v \times \dfrac{C_n^k \times \prod_{i=1}^{T} C_{n-k}^{m_i-k} - \sum_{i=k+1}^{m_{\min}} C_i^k \times C(x_T=i)}{\prod_{i=1}^{T} C_n^{m_i}}$$

$$= \dfrac{\prod_{i=1}^{T} m_i}{n^{T-1}} \times$$

$$\sum_{k=0}^{m_{\min}} ((k-1)+1)^{v-1} \times \dfrac{k \times \dfrac{n}{k} C_{n-1}^{k-1} \times \prod_{i=1}^{T} C_{n-1-(k-1)}^{m_i-1-(k-1)} - \sum_{i=k+1}^{m_{\min}} k \times \dfrac{n}{k} C_{i-1}^{k-1} \times C(x_T=i-1 \mid n-1, M-1)}{\prod_{i=1}^{T} C_{n-1}^{m_i-1}}$$

$$= E(x_T) \sum_{k=0}^{m_{\min}} \sum_{j=0}^{v-1} C_{v-1}^j (k-1)^j \times \dfrac{C_{n-1}^{k-1} \times \prod_{i=0}^{T-1} C_{n-1-(k-1)}^{m_i-1-(k-1)} - \sum_{i=k+1}^{m_{\min}} C_{i-1}^{k-1} \times C(x_T=i-1 \mid n-1, M-1)}{\prod_{i=1}^{T} C_{n-1}^{m_i-1}}$$

$$= E(x_T) \sum_{j=0}^{v-1} C_{v-1}^j \left( \sum_{k=0}^{m_{\min}} (k-1)^j \times \dfrac{C_{n-1}^{k-1} \times \prod_{i=1}^{T} C_{n-1-(k-1)}^{m_i-1-(k-1)} - \sum_{i=k+1}^{m_{\min}} C_{i-1}^{k-1} \times C(x_T=i-1 \mid n-1, M-1)}{\prod_{i=1}^{T} C_{n-1}^{m_i-1}} \right)$$

$$= E(x_T) \sum_{j=0}^{v-1} C_{v-1}^j E(x_T^j \mid n-1, M-1)$$

**Corollary 1**: Any moments (vth) of $p(x_T)$ has the following formula:

$$E(x_T^v) = \sum_{i=1}^{v} g(v,i) \dfrac{\prod_{j=0}^{v-1} \prod_{i=1}^{T} (m_i - j)}{\prod_{j=0}^{v-1} (n-j)^{T-1}}$$

Here $g(v, i)$ are constant coefficients which are only related with $v$ and $i$, and can be determined by the above formulas.

Specifically, for $v=2$, we have:

$$E(x_T^2) = E(x_T)(1+E(x_T \mid n-1, M-1)) = \frac{\prod_{i=1}^{T} m_i}{n^{T-1}} \times \left(1 + \frac{\prod_{i=1}^{T}(m_i-1)}{(n-1)^{T-1}}\right)$$

$$= \frac{\prod_{i=1}^{T} m_i}{n^{T-1}} + \frac{\prod_{i=0}^{T-1} m_i}{n^{T-1}} \times \frac{\prod_{i=1}^{T}(m_i-1)}{(n-1)^{T-1}}$$

For v=3, we have:

$$E(x_T^3) = E(x_T) \times \sum_{i=0}^{2} C_2^i E(x_T^i \mid n-1, M-1)$$

$$= E(x_T) \times \left(1 + 2E(x_T \mid n-1, M-1) + E(x_T^2 \mid n-1, M-1)\right)$$

$$= E(x_T) \times \left(1 + 2E(x_T \mid n-1, M-1) + E(x_T \mid n-1, M-1)(1 + E(x_T \mid n-2, M-2))\right)$$

$$= E(x_T) \times \left(1 + 3E(x_T \mid n-1, M-1) + E(x_T \mid n-1, M-1)E(x_T \mid n-2, M-2)\right)$$

$$= \frac{\prod_{i=1}^{T} m_i}{n^{T-1}} \times \left(1 + \frac{\prod_{i=1}^{T}(m_i-1)}{(n-1)^{T-1}}\left(3 + \frac{\prod_{i=1}^{T}(m_i-2)}{(n-2)^{T-1}}\right)\right)$$

**Corollary 2**: Any central moments ($v$th) about the mean can be calculated by:

$$E\left((x_T - E(x_T))^v\right) = \sum_{j=0}^{v} C_v^j (-1)^{v-j} E(x_T^j)(E(x_T))^{v-j}$$

**Corollary 3**: let $T = 2$, then we can get any moment and central moment for classical HGD.

**Corollary 4**: let $v=2$, we can get the variance for the $p(x_T=k)$:

$$var(x_T) = E(x_T^2) - (E(x_T))^2$$

$$= E(x_T)(1 + E(x_T \mid n-1, M-1) - E(x_T))$$

$$= \frac{\prod_{i=1}^{T} m_i}{n^{T-1}} \times \left(1 + \frac{\prod_{i=1}^{T}(m_i-1)}{(n-1)^{T-1}} - \frac{\prod_{i=1}^{T} m_i}{n^{T-1}}\right)$$

Calculation of higher moments and central moments significantly characterized the intrinsic property of the distribution p(k, T), and is useful in probability estimation by taking advantage of different forms of Chebyshev's inequality.

**4.2  Appendix 2: Formulas for Mathematical expectation of p(x$_t$) and p(x$_{\geq t}$)**

**We will use the same labels defined in the results section. In brief, $d = T-t$, $M_t$ is one of the $t$-subset (that is a subset containing $t$ element) of M, and { $M_t$ } is the set of all possible $t$-subsets of M. $\overline{M_t}$ is the complementary set of $M_t$ in M, that is, $\overline{M_t} = \{n - \overline{m}_{t_1}, n - \overline{m}_{t_2}, ... n - \overline{m}_{t_d}\}$, here $\overline{m}_{t_i}$ ($1 \leq i \leq d$) is the elements that belong to the M but not to the $M_t$, and $M(M_t)$ is the union set of $M_t$ and $\overline{M_t}$.**

**Theorem 3.1 Mathematical expectation for p(x$_t$) is**

$$E(x_t) = \sum_{M_t \subset \{Mt\}} E(x_T \mid n, M(M_t)) = \sum_{M_t \subset \{Mt\}} \frac{\prod_{i=1}^{T} M(M_t)}{n^{T-1}}$$

when all $m_i = m_0$ ($1 \leq i \leq T$), the Mathematical expectation is:

$$E(x_t) = C_T^t \times \frac{(n-m_0)^{T-t} \times m_0^{t}}{n^{T-1}}$$

**Theorem 3.2 Mathematical expectation for p(x$_{\geq t}$) is**

$$E(x_{\geq t}) = E(x_t) + E(x_{\geq (t+1)}) = \sum_{i=t}^{T} E(x_i)$$

**Proof:**

The key point for the formulas is that, the variable $x_t$ is the sum of variables $x_t(M_t)$. Here $x_t(M_t)$ is the number of $t$-overlapped elements in one $t$-subset ($M_t$) of M. In detail, we have:

$$x_t = \sum_{M_t \subset \{M_t\}} x_t(M_t) \quad [A2]$$

Notably, the $x_t(M_t)$ can be viewed as the $T$-overlapped elements in $M(M_t)$, and whose mathematical expectation can be calculated by Theorem 1. Similarly, the variable $x_{\geq t}$ is the sum of $x_t$ and $x_{\geq t+1}$, leading to the result of Theorem 3.2.

***Variance of p(x$_t$) and p(x$_{\geq t}$)***

**Lemma 2**. The mathematical expectation of $Y = E(x_T \times x_0)$ is:

$$E(y = x_T \times x_0) = E(y = x_T \times x_0 \mid n, M) = \frac{\prod_{i=1}^{T} m_i \prod_{i=1}^{T} (n - m_i)}{(n-1)^{T-1} n^{T-1}}$$

**Lemma 3**. Select 2 subsets, *M1*, *M2*, which are not overlapped in *N*, and select another subset *M*, the mathematical expectation of y=($x_{2,MM1} \times x_{2,MM2}$) is: (Here, $x_{2,MM1}$ and $x_{2,MM2}$ is the number of overlapped elements between *M* and *M1* and *M2*, respectively)

$$E(y = x_{2,MM1} \times x_{2,MM2} \mid n,m,m1,m2) = \frac{m1 \times m2 \times m \times (m-1)}{n(n-1)}$$

**Lemma 4.** Select $m_i$ (1⩽i⩽T) from *N*, and M={ $m_1, m_2...m_T$ }, let $M_t$ is a *t*-subset in *M*, and $M_s$ is a *s*-subset in *M*, let $x_{t,Mt}$ and is the *t*-overlapped elements just in the $M_t$ and $x_{s,Ms}$ is the *s*-overlapped elements just in $M_s$, then the mathematical expectation of $y = x_{t,Mt} \times x_{s,Ms}$ is:

$$E(y = x_{t,Mt} \times x_{s,Ms}) = \frac{E(y_{M_{NoSame}} = x_{T',M_{noSame}} \times x_{0,M_{noSame}}) \times (E(x^2_{T'',M\_same}) - E(x_{T'',M\_same}))}{n(n-1)}$$

$$= \frac{\prod_{m_i \subset M_{noSame}} m_i(n-m_i) \prod_{m_j \subset M_{same}} m_j(m_j-1)}{n^T(n-1)^{T-1}}$$

Here, the key of the formulas is the definition of $M_{same}$ and $M_{noSame}$. In detail, $M_{same}$ is the intersection set of $M(M_t)$ and $M(M_s)$, and $M_{noSame}$ is the elements in $M(M_t)$ but not in $M_{same}$.

**Theorem 4** According to equation [A2], and lemma 2~4, we have the 2$^{nd}$ 0-moment of $p(x_t)$ as:

$$E(x_t^2) = \sum_{M_{t1},M_{t2} \subset \{M_t\}} E(x_{t,Mt1} \times x_{t,Mt2})$$

$$= \sum_{\substack{M_{t1},M_{t2} \subset \{M_t\} \\ M_{t1} \neq M_{t2}}} \frac{\prod_{m_i \subset M_{noSameMt1Mt2}} m_i(n-m_i) \prod_{m_j \subset M_{sameMt1Mt2}} m_j(m_j-1)}{n^T(n-1)^{T-1}} + \sum_{\substack{M_{t1},M_{t2} \subset \{M_t\} \\ M_{t1}=M_{t2}}} \frac{\prod_{m_i \subset (M(M_{t1}))} m_i}{n^{T-1}} \left(1 + \frac{\prod_{m_i \subset (M(M_{t1}))}(m_i-1)}{(n-1)^{T-1}}\right)$$

### 4.3 *Appendix 3 Mathematical Proof for E($x_T$) > Var($x_T$) and the difference between E($x_T$) and Var($x_T$) is a higher order infinitesimal of E($x_T$) when $E(x_T) \to 0$.*

We will proof this for the situation fully overlapped elements $p(x_T = k)$. We have:

$$Var(x_T) = \frac{\prod_{i=1}^T m_i}{n^{T-1}} \times \left(1 + \frac{\prod_{i=1}^T(m_i-1)}{(n-1)^{T-1}} - \frac{\prod_{i=1}^T m_i}{n^{T-1}}\right) = E(x_T) \times \left(1 + \frac{\prod_{i=1}^T(m_i-1)}{(n-1)^{T-1}} - \frac{\prod_{i=1}^T m_i}{n^{T-1}}\right)$$

We will proof that, $\left(1 + \frac{\prod_{i=1}^T(m_i-1)}{(n-1)^{T-1}} - \frac{\prod_{i=1}^T m_i}{n^{T-1}}\right) < 1$, then we have $E(x_t) > Var(x_t)$.

Because $M[i] \leq N$, then we have $\frac{\prod_{i=1}^T(m_i-1)}{\prod_{i=1}^T m_i} = \prod_{i=1}^T \left(1 - \frac{1}{m_i}\right) < \prod_{i=1}^{T-1}\left(1 - \frac{1}{m_i}\right) < \prod_{i=1}^{T-1}\left(1 - \frac{1}{n}\right) = \frac{(n-1)^{T-1}}{n^{T-1}}$

That is $\frac{\prod_{i=1}^T(m_i-1)}{\prod_{i=1}^T m_i} < \frac{(n-1)^{T-1}}{n^{T-1}}$, $\frac{\prod_{i=1}^T(m_i-1)}{(n-1)^{T-1}} < \frac{\prod_{i=1}^T m_i}{n^{T-1}}$, then we get the result.

Next, $\frac{E(x_T) - Var(x_T)}{E(x_T)} = \frac{\prod_{i=1}^{T} m_i}{n^{T-1}} - \frac{\prod_{i=1}^{T}(m_i - 1)}{(n-1)^{T-1}} < E(x_T)$, obviously, when $E(x_T) \to 0$, it approaches to 0.

Therefore, the difference between $E(x_T)$ and $Var(x_T)$ is a higher order infinitesimal of $E(x_T)$ when $E(x_T) \to 0$.

## ACKNOWLEDGMENTS

Funding and Conflict of interest: This study was partially supported by National Natural Science Foundation of China (81972359, 81502143).